\documentclass{amsart}
\usepackage[cp866]{inputenc}
\usepackage[english,russian]{babel}

\usepackage{amssymb}
\usepackage{amscd}
\begin{document}
\newcommand\Rs{{\mathbb R}}
\newcommand\End{\mathop {\fam 0 End}\nolimits}
\newcommand\Diff{\mathop {\fam 0 Diff}\nolimits}
\newcommand\id{\mathop {\fam 0 id}\nolimits}
\newcommand\Der{\mathop {\fam 0 Der}\nolimits}
\newcommand\ks{{\Bbbk}}

\title{ФОРМУЛА ЛЕЙБНИЦА ДЛЯ КОВАРИАНТНОЙ ПРОИЗВОДНОЙ
И НЕКОТОРЫЕ ЕЕ ПРИЛОЖЕНИЯ}
\thanks{Работа частично поддержана грантом РФФИ-08-01-92001}
\author{А.~В.~Гаврилов}

\maketitle
\smash{\raise 70mm\hbox{УДК 514.76}}

{\small Получена формула для высших ковариантных производных на тензорном
произведении векторных расслоений, являющаяся  широким обобщением
классической формулы Лейбница. Построен алгоритм вычисления линейной по
второй переменной части ряда Тейлора двойного экспоненциального отображения.
}

{\small Ключевые слова и фразы: аффинная связность, связность на векторном
расслоении, высшие ковариантные производные,
композиция экспоненциальных отображений.}

\section{Введение}

Пусть ${\mathcal M}$ -- гладкое многообразие и $E\to {\mathcal M}$ --
гладкое векторное расслоение. Через
${\mathfrak F}={\mathfrak F}({\mathcal M})$
мы обозначаем алгебру гладких функций на ${\mathcal M}.$
${\mathfrak F}$ -- модуль гладких сечений $E$
мы обозначаем $\Gamma({\mathcal M},E)$. Для модуля сечений касательного
расслоения, т.е. для пространства
векторных полей, используется обозначение
${\mathcal V}={\mathcal V}({\mathcal M})=\Gamma({\mathcal M},T{\mathcal M})$.
Символы $T^n$ и $T_n$ означают соответственно контравариантные и ковариантные
тензорные расслоения соответствующего ранга:
$T^n{\mathcal M}=T{\mathcal M}^{\otimes n},\,
T_n{\mathcal M}=T^*{\mathcal M}^{\otimes n}$.
Через
$\Diff(\Gamma({\mathcal M},E))$
будем обозначать алгебру дифференциальных операторов, действующих
в пространстве сечений $E.$

Ковариантной производной на векторном расслоении $E$ называется линейный
оператор
$$\nabla:\Gamma({\mathcal M},E)\to
\Gamma({\mathcal M},T_1{\mathcal M}\otimes E),$$
для которого имеет место тождество
$$\nabla(fu)=f\nabla(u)+df\otimes u,\,f\in{\mathfrak F},
u\in\Gamma({\mathcal M},E).$$
Связность на векторном расслоении -- то же, что и ковариантная
производная на нём (для связности на более
общих расслоениях такое определение непригодно).
Аффинную связность на многообразии ${\mathcal M}$ можно рассматривать
как частный случай, а именно как связность на
кокасательном расслоении $T_1{\mathcal M}$ (или, что то же самое,
на касательном расслоении). Если имеются два расслоения со связностями
$E_1, E_2\to {\mathcal M}$, то на их тензорном произведении $E_1\otimes E_2$
естественным образом вводится связность
$$\nabla:\Gamma({\mathcal M},E_1\otimes E_2)\to
\Gamma({\mathcal M},T_1{\mathcal M}\otimes E_1\otimes E_2),$$
согласно правилу
$$\nabla u_1\otimes u_2=(\nabla u_1)\otimes u_2+u_1\otimes(\nabla u_2).$$

Предположим, что на ${\mathcal M}$ имеется аффинная связность и расслоение
$E\to {\mathcal M}$
также снабжено связностью. Тогда оператор ковариантной производной
распространяется на все расслоения вида $T_n{\mathcal M}\otimes E$:
$$\nabla:\Gamma({\mathcal M},T_n{\mathcal M}\otimes E)\to
\Gamma({\mathcal M},T_{n+1}{\mathcal M}\otimes E),\,n\ge 0.$$
По определению, высшими ковариантными производными называются степени этого
оператора:
$$\nabla^n:\Gamma({\mathcal M},E)\to
\Gamma({\mathcal M},T_{n}{\mathcal M}\otimes E).$$

Классическая формула Лейбница выражает высшие производные произведения через
производные сомножителей:
$$(fg)^{(n)}=\sum_{m=0}^n\binom{n}{m}f^{(m)}g^{(n-m)},\,
f,g\in{\mathfrak F}(\Rs^d).$$
Одной из целей настоящей работы является обобщение
формулы Лейбница на ковариантные производные на расслоениях.
Именно, нами будет получено выражение производной $\nabla^n (u_1\otimes u_2),$
где $u_1\in\Gamma({\mathcal M},E_1),\,u_2\in\Gamma({\mathcal M},E_2),$
через производные вида $\nabla^m u_1\otimes\nabla^{n-m} u_2$.
В качестве следствия будет получено новое доказательство теоремы 1 [3],
более геометричное по своему характеру и более прозрачное в сравнении
с первоначальным.

Во второй части работы рассматривается приложение этой формулы к задаче,
связанной с двойным экспоненциальным отображением. Речь идет о вычислении
ряда Тейлора введённой в [1] функции $h(v,w)$ по переменным $v,w$.
Судя по всему, в целом это весьма сложная проблема. Однако, используя
алгебраические методы наряду с геометрическими соображениями,
удается построить алгоритм, позволяющий находить линейные по второй переменной
слагаемые ряда с достаточной эффективностью.

\section{Ковариантная формула Лейбница}

Хотя ковариантное обобщение формулы Лейбница не очень
сложно по своему содержанию,
его формулировка не так уж очевидна. С формальной точки зрения,
смысл этой формулы состоит в
явном описании определённых сечений расслоений вида
$${\mathop {\fam 0 Hom}\nolimits}
(T_m{\mathcal M}\otimes E_1\otimes T_{n-m}{\mathcal M}\otimes E_2,
T_n{\mathcal M}\otimes E_1\otimes E_2).$$
Довольно трудно подобрать подходящий для такого описания язык.

Простейший, по мнению автора, путь
заключается в использовании высших производных в сочетании со сверткой
с контравариантным тензором.
Если $A\in\Gamma({\mathcal M},T^n{\mathcal M})$ есть тензорное поле
ранга $n,$ то определена свертка $A$ с $\nabla^n,$ являющаяся оператором
вида
$$A\cdot\nabla^n:u\mapsto A\cdot(\nabla^n u),\,u\in\Gamma({\mathcal M},E).$$
Применяя эту конструкцию на практике, необходимо следить за порядком инексов в
$A$ и $\nabla^n$ во избежание путаницы ( этот вопрос подробно
рассматривается в книге Пале, указанной в библиографии к [3]).
Для свертки с тензорным произведением
векторных полей обычно используется обозначение
$$\nabla^n_{v_1,\dots,v_n}=(v_1\otimes\dots\otimes v_n)\cdot\nabla^n,\,
v_1,\dots,v_n\in{\mathcal V}.$$

Следуя [3], мы будем обозначать оператор $A\cdot\nabla^n$,
действующий в расслоении $E$, через $\mu_E(A)\in\Diff(\Gamma({\mathcal M},E))$.
Как и в [3], мы будем использовать обозначение
$$T({\mathcal M})=\bigoplus_{n\ge 0}\Gamma({\mathcal M},T^n{\mathcal M})$$
для формальной суммы ${\mathfrak F}$ - модулей сечений тензорных расслоений
(считая $\Gamma({\mathcal M},T^0{\mathcal M})={\mathfrak F}$).
Элементами $T({\mathcal M})$ являются
(конечные) формальные суммы тензорных полей различного ранга.
Этот модуль является ${\mathfrak F}$ - алгеброй
относительно тензорного умножения.

Нам будет полезна следующая

{\bf Лемма 1}
{\it Пусть на многообразии имеются два векторных расслоения

$E_1, E_2\to {\mathcal M}$. Тогда естественное отображение
$$\Gamma({\mathcal M},E_1)\underset{{\mathfrak F}}{\otimes}
\Gamma({\mathcal M},E_2)\to\Gamma({\mathcal M},E_1\otimes E_2)$$
является изоморфизмом ${\mathfrak F}$ - модулей.
}

Утверждение леммы хорошо известно и легко доказывается, например, при
помощи разбиения единицы. Однако автор, к сожалению,
не нашел нужной формулировки
в литературе. В случае, когда $E_1$ и $E_2$ - ковариантные тензорные
расслоения, оно эквивалентно предложению 3.1 из [4, гл I].
Одним из следствий леммы является канонический изоморфизм
$$T({\mathcal M})\cong T_{\mathfrak F}({\mathcal V}),$$
который будет использован далее.

Оператор $\mu_E$ на расслоении со связностью можно рассматривать как
гомоморфизм ${\mathfrak F}$ - модулей:
$$\mu_E:T({\mathcal M})\to\Diff(\Gamma({\mathcal M},E)),\,
\mu_E:A\mapsto A\cdot\nabla^n.$$
Следует помнить, что он {\it не} является гомоморфизмом алгебр.
Приводимая ниже формула (4) как раз показывает, в какой мере $\mu_E$
отличается от гомоморфизма алгебр (над полем $\Rs$).

Если имеются два расслоения со связностью $E_1, E_2\to {\mathcal M}$,
то каждому элементу $A\in T({\mathcal M})$ можно поставить в соответствие
три оператора, действующие в разных расслоениях:
$$\mu_{E_1}(A)\in\Diff(\Gamma({\mathcal M},E_1)),\,
\mu_{E_2}(A)\in\Diff(\Gamma({\mathcal M},E_2)),\,
\mu_{E_1\otimes E_2}(A)\in\Diff(\Gamma({\mathcal M},E_1\otimes E_2)).$$

Ковариантная формула Лейбница устанавливает связь между этими тремя операторами.
Чтобы сформулировать её, нам потребуется введенное в [3]
коумножение в алгебре $T({\mathcal M})$, которое является
${\mathfrak F}$ -  линейным отображением вида
$$\bigtriangleup:T({\mathcal M})\to
T({\mathcal M})
\underset{{\mathfrak F}}{\breve{\otimes}}
T({\mathcal M}).$$
Как и в [3], мы используем знак $\breve{\otimes}$ для того,
чтобы отличать соответствующее умножение от умножения в
алгебре $T({\mathcal M})$ (обозначаемого, как обычно, символом $\otimes$).
Кроме того, мы используем
обозначение ${\mathfrak t}_{E_1\otimes E_2}$ для естественного
отображения перехода к тензорному произведению над ${\mathfrak F}$:
$${\mathfrak t}_{E_1\otimes E_2}:
\Gamma({\mathcal M},E_1)\underset{\Rs}{\otimes}\Gamma({\mathcal M},E_2)\to
\Gamma({\mathcal M},E_1)\underset{{\mathfrak F}}{\otimes}
\Gamma({\mathcal M},E_2)=\Gamma({\mathcal M},E_1\otimes E_2).$$

{\bf Теорема 1}
{\it Для всякого $A\in T({\mathcal M})$ следующая диаграмма
коммутативна:
\[
\begin{CD}
\Gamma({\mathcal M},E_1)
\underset{\Rs}{\otimes}\Gamma({\mathcal M},E_2)
@>\sum_{(A)}\mu_{E_1}(A_{(1)})\otimes\mu_{E_2}(A_{(2)})>>
\Gamma({\mathcal M},E_1)
\underset{{\mathfrak F}}{\otimes}\Gamma({\mathcal M},E_2)\\
@VV {\mathfrak t}_{E_1\otimes E_2} V    @VV  V \\
\Gamma({\mathcal M},E_1\otimes E_2)
@>\mu_{E_1\otimes E_2}(A)>>
\Gamma({\mathcal M},E_1\otimes E_2)
\end{CD}
\]
}

Здесь использованы обычные обозначения Свидлера.
Отображение в верхней строке диаграммы действует следующим образом:
$$\sum_{(A)}\mu_{E_1}(A_{(1)})\otimes\mu_{E_2}(A_{(2)}):
u_1\underset{\Rs}{\otimes}u_2\mapsto
\sum_{(A)}\mu_{E_1}(A_{(1)})u_1
\underset{{\mathfrak F}}{\otimes}\mu_{E_2}(u_{(2)})u_2,$$
где $u_1\in\Gamma({\mathcal M},E_1),\,u_2\in\Gamma({\mathcal M},E_2)$ и
$\bigtriangleup(A)=\sum_{(A)}A_{(1)}\breve{\otimes}A_{(2)}$.
Причина, по которой диаграмма начинается с тензорного произведения
${\mathfrak F}$ - модулей над вещественным полем, заключается в том, что
операторы
$\mu_{E_1}(A_{(1)})$ и $\mu_{E_2}(A_{(2)})$
не являются  ${\mathfrak F}$ - линейными, поэтому отображение
в верхней строке нельзя, по формальным причинам, применять к произведению
$\Gamma({\mathcal M},E_1)
\underset{{\mathfrak F}}{\otimes}\Gamma({\mathcal M},E_2)$.

Утверждение теоремы можно записать в форме равенства
$$\mu_{E_1\otimes E_2}(A)u_1\otimes u_2=\sum_{(A)}
\mu_{E_1}(A_{(1)})u_1\otimes\mu_{E_2}(A_{(2)})u_2,\,A\in T({\mathcal M}),
\eqno{(1)}$$
которое и является ковариантным обобщение формулы Лейбница.
Чтобы в этом убедиться, рассмотрим случай $A=v^{\otimes n},$
где $v\in{\mathcal V}$ - векторное поле. Так как
$\bigtriangleup(v^{\otimes n})=\sum_{m=0}^{n}\binom{n}{m}
v^{\otimes m}\breve{\otimes}v^{\otimes n-m}$ [3, \S3],
мы получаем
$$\mu_{E_1\otimes E_2}(v^{\otimes n})u_1\otimes u_2=\sum_{m=0}^n\binom{n}{m}
\mu_{E_1}(v^{\otimes m})u_1\otimes\mu_{E_2}(v^{\otimes n-m})u_2.$$
По линейности последнее равенство продолжается на симметричные тензоры
$A\in T({\mathcal M})$; следуя [3], мы можем переписать его
в мультииндексной форме (в локальных координатах):
$$\nabla^{\alpha}u_1\otimes u_2=\sum_{\beta}\binom{\alpha}{\beta}
\nabla^{\beta}u_1\otimes\nabla^{\alpha-\beta}u_2.\eqno{(2)}$$
Если оба расслоения тривиальны, т.е. $E_1=E_2=\Rs\times{\mathcal M}$,
то их сечения -- обычные функции, и (2) превращается в знакомую формулу
Лейбница (в мультииндексной записи):
$$\nabla^{\alpha}(fg)=\sum_{\beta}\binom{\alpha}{\beta}
(\nabla^{\beta}f)(\nabla^{\alpha-\beta}g),\,
f,g\in{\mathfrak F}.$$

\section{Доказательство формулы Лейбница}

Наш подход к доказательству ковариантной формулы Лейбница связан
с использованием специального отображения, называемого
$K$ - оператором\footnote{Это не очень удачное название автор произвёл от
слова "ковариантный" (так что буква "$K$"\-  из русского алфавита).
Весьма возможно, что сама конструкция не нова, но автор не сумел найти её
в литературе}. Напомним его определение [3].
Пусть ${\mathfrak g}$  -- неассоциативная алгебра над полем $\ks$
характеристики нуль; следуя [1,3] будем обозначать операцию в алгебре
символом $\diamond.$
Пусть $T({\mathfrak g})$ - тензорная алгебра ${\mathfrak g}$ как линейного
пространства.
Для $x\in{\mathfrak g}$ будем обозначать
через $\tau_x$ отображение вида $\tau_x:y\mapsto x\diamond y$, продолженное
до дифференцирования $T({\mathfrak g})$.

Мы называем $K$ - оператором отображение
$K:T({\mathfrak g})\to T({\mathfrak g})$,
однозначно определяемое условиями
$K(1)=1$ и
$$K(x\otimes a+\tau_x a)=x\otimes K(a),\,x\in{\mathfrak g},
a\in T({\mathfrak g}).$$ Нетрудно убедиться, что это отображение обратимо.
Как известно, тензорная алгебра $T({\mathfrak g})$ обладает структурой
алгебры Хопфа, со стандартным коумножением
$$\bigtriangleup:T({\mathfrak g})\to
T({\mathfrak g})\breve{\otimes}T({\mathfrak g}),\,
\bigtriangleup:x\mapsto 1\breve{\otimes}x+x\breve{\otimes}1,\,
x\in{\mathfrak g}$$
(знак $\breve{\otimes}$ используется из тех же соображений, что и раньше).

{\bf Предложение 1}
{\it Для всякой неассоциативной алгебры ${\mathfrak g}$ диаграмма
\[
\begin{CD}
T({\mathfrak g}) @> K >> T({\mathfrak g})  \\
@VV \bigtriangleup V	@VV \bigtriangleup V \\
T({\mathfrak g})\breve{\otimes}T({\mathfrak g})
@> K\breve{\otimes}K  >>
T({\mathfrak g})\breve{\otimes}T({\mathfrak g})
\end{CD}
\]
коммутативна.
}

Другими словами, $K$ является автоморфизмом $T({\mathfrak g})$ как коалгебры.
В обозначениях Свидлера утверждение записывается в виде равенства
$$\bigtriangleup(K(a))=\sum_{(Ka)}
(Ka)_{(1)}\breve{\otimes}(Ka)_{(2)}=
\sum_{(a)}K(a_{(1)})\breve{\otimes}K(a_{(2)})=
K\breve{\otimes}K(\bigtriangleup(a)),\eqno{(3)}$$
где $a\in T({\mathfrak g})$.
Например, для $x,y\in{\mathfrak g}$ мы имеем
$$K\breve{\otimes}K(\bigtriangleup(x\otimes y))=
K\breve{\otimes}K(1\breve{\otimes}x\otimes y+
x\breve{\otimes}y+y\breve{\otimes}x+x\otimes y\breve{\otimes}1)=
1\breve{\otimes}(x\otimes y-x\diamond y)+
x\breve{\otimes}y+$$
$$+y\breve{\otimes}x+
(x\otimes y-x\diamond y)\breve{\otimes}1=
\bigtriangleup(K(x\otimes y)),$$
поскольку $K(x)=x$ и $K(x\otimes y)=x\otimes y-x\diamond y$.

Для доказательства нам потребуется следующая лемма,
без труда доказываемая по индукции:

{\bf Лемма 2}
{\it  Пусть $V$ -- линейное пространство над полем $\ks$,
$\delta\in\Der_{\ks}(T(V))$ и
$a\in T(V)$. Тогда
$$\bigtriangleup(\delta a)=
\sum_{(a)}(\delta a_{(1)})\breve{\otimes}a_{(2)}+
a_{(1)}\breve{\otimes}(\delta a_{(2)}),$$
где $\bigtriangleup(a)=\sum_{(a)}a_{(1)}\breve{\otimes}a_{(2)}.$
}

Равенство (3) мы будем доказывать индукцией по степени $a$. Для
$a\in{\mathfrak g}$ оно очевидно. Предположим, что равенство справедливо
при условии $\deg(a)\le n$. Пусть $a$ -- элемент $T({\mathfrak g})$,
удовлетворяющий этому условию, и $x\in{\mathfrak g}$.
Обозначим $b=x\otimes a$. Тогда
$$\bigtriangleup(K(b))=
\bigtriangleup(x\otimes K(a)-K(\tau_x a))=
\bigtriangleup(x)\otimes\bigtriangleup(K(a))-\bigtriangleup(K(\tau_x a)).$$
Поскольку $\deg\tau_x a\le \deg a\le n,$  по предположению индукции мы
имеем (опуская знак суммы)
$$\bigtriangleup(K(b))=(1\breve{\otimes}x+x\breve{\otimes}1)\otimes
(K(a_{(1)})\breve{\otimes}K(a_{(2)}))-
K\breve{\otimes}K(\bigtriangleup(\tau_x a)).$$
Применяя лемму 2 к дифференцированию $\tau_x$, мы получим
$$\bigtriangleup(K(b))=x\otimes K(a_{(1)})\breve{\otimes}K(a_{(2)})+
K(a_{(1)})\breve{\otimes}x\otimes K(a_{(2)})-
K(\tau_x a_{(1)})\breve{\otimes}K(a_{(2)})-$$
$$-K(a_{(1)})\breve{\otimes}K(\tau_x a_{(2)})=
K(x\otimes a_{(1)})\breve{\otimes}K(a_{(2)})+
K(a_{(1)})\breve{\otimes}K(x\otimes a_{(2)})=
K\breve{\otimes}K(\bigtriangleup(b)).$$
Таким образом, формула (3) оказывается справедливой
для $b=x\otimes a$. По линейности она распространяется на все
многочлены $b,\,\deg(b)\le n+1$. Следовательно,
формула справедлива для многочленов любой степени.

Перейдём к геометрии. Пусть $T({\mathcal V})$ --
(формальная) тензорная алгебра ${\mathcal V}={\mathcal V}({\mathcal M})$
как вещественного линейного пространства.
Для всякого расслоения $E\to {\mathcal M}$ однозначно определён гомоморфизм
алгебр (над вещественным полем)
$\tau_E:T({\mathcal V})\to\Diff(\Gamma({\mathcal M},E))$
такой, что $\tau_E:v\mapsto\nabla_v$
(к примеру, $\tau_E:v_1\otimes v_2\mapsto\nabla_{v_1}\nabla_{v_2}$).

{\bf Лемма 3}
{\it Пусть имеются два расслоения со связностями $E_1, E_2\to {\mathcal M}$.
Тогда для всякого $a\in T({\mathcal V})$ следующая диаграмма коммутативна:
\[
\begin{CD}
\Gamma({\mathcal M},E_1)\underset{\Rs}{\otimes}\Gamma({\mathcal M},E_2)
@> \sum_{(a)}\tau_{E_1}(a_{(1)})\otimes\tau_{E_2}(a_{(2)})>>
\Gamma({\mathcal M},E_1)\underset{\Rs}{\otimes}\Gamma({\mathcal M},E_2)\\
@VV {\mathfrak t}_{E_1\otimes E_2} V @VV {\mathfrak t}_{E_1\otimes E_2} V \\
\Gamma({\mathcal M},E_1\otimes E_2)
@>\tau_{E_1\otimes E_2}(a)>>
\Gamma({\mathcal M},E_1\otimes E_2)
\end{CD}
\]
}

Здесь $\bigtriangleup(a)=\sum_{(a)}a_{(1)}\breve{\otimes}a_{(2)}$
обозначает стандартное коумножение в тензорной алгебре.
Доказательство леммы несложно. В случае $a=1$ утверждение очевидно, поскольку
$\tau_E(1)$ - тождественный оператор. Для $a=v\in{\mathcal V}$
оно является следствием равенства
$$\nabla_v(u_1\otimes u_2)=(\nabla_v u_1)\otimes u_2+u_1\otimes (\nabla_v u_2),$$
вытекающего из определения ковариантной производной на произведении
расслоений.

Заметим, что как отображения $\tau_E$, так и коумножение являются
гомоморфизмами алгебр. Отсюда, как нетрудно видеть, следует, что если диаграмма
коммутативна при $a=b_1$ и при $a=b_2$, то она будет коммутативной и
для $a=b_1\otimes b_2$.
По индукции мы получаем, что она коммутативна для любого $a$.

Теперь мы можем перейти к доказательству теоремы.
Рассмотрим пространство векторных полей ${\mathcal V}$ как
неассоциативную алгебру с операцией, задаваемой ковариантной производной:
$v\diamond w=\nabla_vw,\,v,w\in{\mathcal V}$.
Тогда на тензорной алгебре $T({\mathcal V})$ определён соответствующий
$K$ - оператор. Он связан с высшими ковариантными производными на расслоении
равенством $\mu_E\circ{\mathfrak t}=\tau_E\circ K$ [3, лемма 2],
где
$${\mathfrak t}:T({\mathcal V})\to
T_{{\mathfrak F}}({\mathcal V})=T({\mathcal M})$$
есть операция перехода
к тензорному умножению над ${\mathfrak F}$.
(Заметим, что пространство ${\mathcal V}$ наделено тремя структурами:
алгебры Ли, ${\mathfrak F}$ - модуля и неассоциативной алгебры.
В данной статье структура алгебры Ли игнорируется. В других
работах автора [1,6] она, напротив, играет главную роль).

Пусть $A={\mathfrak t}(a),\,a\in T(V)$. По определению коумножения в
$T({\mathcal M})$,
$$\bigtriangleup(A)=\sum_{(a)}{\mathfrak t}(a_{(1)})
\underset{{\mathfrak F}}{\breve{\otimes}}{\mathfrak t}(a_{(2)}).$$
Следовательно, диаграмму из условия теоремы можно переписать в виде

\[
\begin{CD}
\Gamma({\mathcal M},E_1)
\underset{\Rs}{\otimes}\Gamma({\mathcal M},E_2)
@>\sum_{(a)}\tau_{E_1}(K(a_{(1)}))\otimes\tau_{E_2}(K(a_{(2)}))>>
\Gamma({\mathcal M},E_1)
\underset{{\mathfrak F}}{\otimes}\Gamma({\mathcal M},E_2)\\
@VV {\mathfrak t}_{E_1\otimes E_2} V    @VV  V \\
\Gamma({\mathcal M},E_1\otimes E_2)
@>\tau_{E_1\otimes E_2}(K(a))>>
\Gamma({\mathcal M},E_1\otimes E_2)
\end{CD}
\]

Так как, согласно предложению 1,
$$\bigtriangleup(K(a))=
\sum_{(a)}K(a_{(1)})\breve{\otimes}K(a_{(2)}),$$
диаграмма коммутативна по лемме 3.

\section{Тождество коцикла}

Одним из основных результатов [3] была формула
$$\mu_E(A)\mu_E(B)=\sum_{(A)}\mu_E(A_{(1)}\otimes\mu_T(A_{(2)})B)\eqno{(4)}$$
[3, теорема 1].
Ее доказательство, основанное на рассматриваемом ниже алгебраическом тождестве,
было довольно непрозрачным. Пользуясь ковариантной
формулой Лейбница, можно дать другое доказательство этого равенства, более
геометрического характера.
Пусть
$A\in\Gamma({\mathcal M},T^n{\mathcal M}),\,
B\in\Gamma({\mathcal M},T^m{\mathcal M})$.
Тогда
$$\mu_E(B)u=B\cdot\nabla^mu=C(B\otimes\nabla^mu),\,
u\in\Gamma({\mathcal M},E),$$
где символом $C$ обозначена свёртка. Применяя (1) к расслоениям
$E_1=T^m{\mathcal M},\,E_2=T_m{\mathcal M}\otimes E,$ мы получим
$$\mu_{E_1\otimes E_2}(A)(B\otimes\nabla^mu)=\sum_{(A)}
\mu_{E_1}(A_{(2)})B\otimes\mu_{E_2}(A_{(1)})\nabla^mu=
\sum_{(A)}\mu_T(A_{(2)})B\otimes[A_{(1)}\cdot\nabla^{\deg(A_{(1)})}
\nabla^m u].$$
Здесь удобно переставить тензоры $A_{(1)}$ и $A_{(2)}$,
пользуясь кокоммутативностью.
Можно, не нарушая общности, считать эти тензоры однородными;
при этом $\deg(A_{(1)})+\deg(A_{(2)})=\deg(A)=n$.

Как известно, свёртка перестановочна с ковариантной производной.
Следовательно,
$$\mu_E(A)\mu_E(B)u=\mu_E(A)C(B\otimes\nabla^mu)=
C\mu_{E_1\otimes E_2}(A)(B\otimes\nabla^mu)=$$
$$=C\sum_{(A)}\mu_T(A_{(2)})B\otimes[A_{(1)}\cdot
\nabla^{\deg(A_{(1)}})\nabla^mu]=
\sum_{(A)}[A_{(1)}\otimes\mu_T(A_{(2)})B]\cdot
\nabla^{\deg(A_{(1)})}\nabla^mu=$$
$$=\sum_{(A)}\mu_E(A_{(1)}\otimes\mu_T(A_{(2)})B)u.$$
Таким образом, формула доказана.
Заметим, что при действии свертки
тензоры $A_{(1)}$ и $\mu(A_{(2)})B$ поменялись местами.
Это связано с правилами построения "свёрнутых" ковариантных
производных: тензор $A_{(1)}$ сворачивается с левыми индексами
производной $\nabla^{\deg(A_{(1)})+m}u$.

Как уже упоминалось, в статье [3] формула (4) выводилась из
алгебраического тождества [3, лемма 1]:
$$K(A)\otimes K(B)=\sum_{(A)}K(A_{(1)}\otimes\tau(K(A_{(2)}))B),\eqno{(5)}$$
где через $\tau:T({\mathfrak g})\to\End_{\ks}(T({\mathfrak g}))$
обозначается гомоморфизм
алгебр, определяемый условием $\tau(x)=\tau_x,\,(x\in{\mathfrak g})$
\footnote{Обозначения в настоящей статье слегка отличаются
от использовавшихся в [3]: мы пишем $\tau$ вместо $\hat{\tau}$ и
$\tau_E$ вместо $\tau$. Использование для $\tau$ и $\tau_E$ сходных обозначений
может быть оправдано тем, что между этими двумя операциями
имеется тесная связь [3, лемма 3 ].}.
По причинам, объясняемым ниже, мы будем называть (5)
{\it тождеством коцикла}.

Равенство (5) выгдядит довольно необычно.
Чтобы понять его смысл,
необходимо прежде всего разобраться, что представляет собой гомоморфизм $\tau$
с точки зрения алгебры. Эта операция обладает рядом интересных свойств.
Чтобы их сформулировать, удобно продолжить $\tau$ до гомоморфизма
"удвоенной" тензорной алгебры:
$$\tau:T({\mathfrak g})\breve{\otimes}T({\mathfrak g})\to
\End_{\ks}(T({\mathfrak g})\breve{\otimes}T({\mathfrak g})),\,
\tau(a_1\breve{\otimes}a_2):b_1\breve{\otimes}b_2\mapsto
\tau(a_1)b_1\breve{\otimes}\tau(a_2)b_2,\,a_1,a_2,b_1,b_2\in T({\mathfrak g}).$$

{\bf Предложение 2}
{\it Для всякой неассоциативной алгебры ${\mathfrak g}$
и всякого элемента $a\in T({\mathfrak g})$ следующая диаграмма коммутативна:
\[
\begin{CD}
T({\mathfrak g}) @> \bigtriangleup >>
T({\mathfrak g})\breve{\otimes}T({\mathfrak g}) @> m >>
T({\mathfrak g})  \\
@VV \tau(a) V    @VV \tau(\bigtriangleup(a)) V  @VV \tau(a) V \\
T({\mathfrak g}) @> \bigtriangleup >>
T({\mathfrak g})\breve{\otimes}T({\mathfrak g}) @> m >>
T({\mathfrak g})  \\
\end{CD}
\]
}

Здесь символом $m$ обозначено обычное умножение в тензорной алгебре, т.е.
$m:b_1\breve{\otimes}b_2\mapsto b_1\otimes b_2.$
Коммутативность левого квадрата означает равенство
$$\bigtriangleup(\tau(a)b)=\sum_{(a)}\sum_{(b)}
\tau(a_1)b_1\breve{\otimes}\tau(a_2)b_2.$$
При $a=x\in{\mathfrak g}$ оно следует из леммы 2, поскольку
$\tau(x)=\tau_x$ является дифференцированием и 
$\bigtriangleup(x)=1\breve{\otimes}x+x\breve{\otimes}1$.
В общем случае равенство легко доказывается индукцией по степени $a$.
Коммутативность правого квадрата, т.е. равенство
$$\tau(a)b_1\otimes b_2=\sum_{(a)}
\tau(a_1)b_1\otimes\tau(a_2)b_2$$
может быть доказана аналогично лемме 3.

Напомним некоторые определения из алгебры. Пусть определено действие
группы $G$ на группу $X$, т.е. гомоморфизм групп
$G\to{\mathop {\fam 0 Aut}\nolimits}(X)$. Действие элемента $g\in G$ мы
обозначаем точкой: $g:x\mapsto g\cdot x=g(x),x\in X$. По определению,
$g\cdot(xy)=(g\cdot x)(g\cdot y), g_1\cdot(g_2\cdot x)=(g_1g_2)\cdot x$.
Отображение $\phi:G\to X$ называется {\it коциклом}, если
$$\phi(gh)=\phi(g)(g\cdot\phi(h)),\,g,h\in G.$$

Пусть $\ks$ -- некоторое поле.
Очевидно, всякий коцикл $\phi$ продолжается до линейного отображения
$\phi:\ks G\to\ks X$. Нетрудно видеть, что линейное
отображение такого вида является продолжением коцикла если оно, во-первых,
является гомоморфизмом коалгебр и, во-вторых, удовлетворяет тождеству

$$\phi(ab)=\sum_{(a)}\phi(a_{(1)})(a_{(2)}\cdot\phi(b)),\,a,b\in \ks G.\eqno{(6)}$$

Таким образом, можно называть коциклом отображение одной
биалгебры в другую $\phi:H\to R$, если оно является гомоморфизмом коалгебр
и удовлетворяет равенству (6). Правда, само это равенство имеет смысл
лишь при условии, что определено действие $H$ на $R$. Рассматривая
случай $H=\ks G, R=\ks X$ как определяющий пример, можно определить
действие $H$ на $R$ как линейное отображение $H\otimes R\to R$,
обладающее следующими свойствами:

$$(h_1h_2)\cdot x=h_1\cdot(h_2\cdot x),\,
h\cdot(xy)=\sum_{(h)}(h_{(1)}\cdot x)(h_{(2)}\cdot y),\,
\bigtriangleup(h\cdot x)=\sum_{(h)}\sum_{(x)}
h_{(1)}\cdot x_{(1)}\otimes h_{(2)}\cdot x_{(2)}.$$

Здесь $x,y\in R, h,h_1,h_2\in H$; как и выше, действие обозначено
точкой: $h\otimes x\mapsto h\cdot x$ (точности ради, следует еще добавить
условие $h\cdot 1=\epsilon(h)$).
Первое равенство означает, что $R$ является (левым) $H$ - модулем.
Первые два условия вместе означают, что биалгебра $H$ действует на
алгебру $R$, т.е. что $R$ является $H$ -- модульной алгеброй [7, с. 40].
Действие одной биалгебры на другую назовем {\it согласовынным},
если выполнено третье равенство. Как можно убедиться, в этом случае
смэш-произведение $R\# H$ [7, с. 41] само является биалгеброй.

Предложение 2 означает в точности, что отображение
$T({\mathfrak g})\otimes T({\mathfrak g})\to T({\mathfrak g})$
вида $a\otimes b\mapsto a\cdot b=\tau(a)b$
является согласованным действием биалгебры
$T({\mathfrak g})$ на себя. Обозначив
$a=K(A),\,b=K(B)$ и принимая во внимание равенство
$\bigtriangleup(a)=\sum_{(A)}K(A_{(1)})\otimes K(A_{(2)})$,
мы можем представить (5) в виде
$$K^{-1}(ab)=\sum_{(a)}K^{-1}(a_{(1)})a_{(2)}\cdot K^{-1}(b),\,a,b\in
T({\mathfrak g}).$$
Таким образом, отображение $K^{-1}:T({\mathfrak g})\to T({\mathfrak g})$
является коциклом. Более того, можно убедиться,
что для всякой неассоциативной алгебры
${\mathfrak g}$ существует
единственное  согласованное действие $T({\mathfrak g})$ на себя такое, что
$x\cdot y=x\diamond y$ для $x,y\in{\mathfrak g}$ и единственный коцикл,
соответствующий этому действию.
Таким образом, рассматриваемая конструкция достаточно естественна.

\section{Ковариантная симметризация}

Обозначим через $S({\mathcal M})\subset T({\mathcal M})$
пространство симметричных контравариантных тензорны полей.
Определить его формально можно следующим образом. Пусть
$\imath:S({\mathcal V})\hookrightarrow T({\mathcal V})$ обозначает естественное
вложение симметрической алгебры вещественного пространства
${\mathcal V}$ в тензорную алгебру, т.е. такое линейное отображение,
что $\imath:v^n\mapsto v^{\otimes n} (v\in{\mathcal V},\,n>0)$.
Тогда $S({\mathcal M})={\mathfrak t}\circ\imath (S({\mathcal V}))$.

Рассмотрим отображение
$$\mu:T({\mathcal M})\to\Diff({\mathfrak F}({\mathcal M})),$$
соответствующее тривиальному расслоению $\Rs\times{\mathcal M}$.
Оно переводит тензорные поля в скалярные дифференциальные операторы.
Несложные рассуждения показывают, что ограничение $\mu$ на
$S({\mathcal M})$ является изоморфизмом ${\mathfrak F}$ - модулей
[5, введение]. Обратное отображение мы обозначим через
$$\sigma:\Diff({\mathfrak F}({\mathcal M}))\to S({\mathcal M}),\,
\mu\circ\sigma=\id.$$

Следуя Шарафутдинову [5], назовем поле $\sigma(A)$
{\it геометрическим символом} скалярного
оператора $A\in\Diff({\mathfrak F})$. Следует, правда, заметить, что
в [5] геометрическим символом (или просто символом)
называется не поле $\sigma(A)$, а функция на кокасательном раслоении,
получаемая из него путем свёртки с подходящими степенями ковектора
$\xi$. Автор предлагает различать символ (функцию на кокасательном
расслоении) и геометрический символ. Например, символом
оператора $\nabla^2_{v,v}\, (v\in{\mathcal V})$ будет функция
$\langle v,\xi\rangle^2$, а его геометрический символ
равен $\sigma(\nabla^2_{v,v})=v\otimes v$.

Оператором {\it ковариантной симметризации} будем называть отображение
$$\Theta=\sigma\circ\mu:T({\mathcal M})\to S({\mathcal M}).$$
Из определения следует, что
это проектор на подпространство $S({\mathcal M})$, т.е.
$\Theta^2=\Theta$ и  $\Theta A=A$ при $A\in S({\mathcal M})$.
По мнению автора, явное описание этого проектора было бы полезно
для многих приложений. К сожалению, это, по-видимому,
чрезвычайно сложная задача. В настоящей работе
рассматривается весьма частный случай: симметризация мономов вида
$v^{\otimes n}\otimes w\otimes v^{\otimes m}\,(v,w\in{\mathcal V})$;
но и для этого случая речь идет не о явных формулах, а лишь об
алгоритме, более или менее экономичном.
Геометрический смысл такого рода вычислений объясняется в \S 7.
К этому можно добавить, что поставленная в [5] проблема вычисления
многочленов $R^{\alpha,\beta}$ может быть сведена к вычислению
тензоров $\Theta v^{\otimes n}\otimes w^{\otimes m}$, где
$n=|\alpha|,m=|\beta|$.

Следует сказать, что понятие геометрического символа распространяется и на
операторы, действующие на векторном расслоении [5, II]. Ради полноты,
рассмотрим эту конструкцию. При $\dim E>1$ отображение
$$\mu_E:T({\mathcal M})\to\Diff(\Gamma({\mathcal M},E))$$
не сюрьективно. Чтобы выйти из положения, нужно ввести в рассмотрение модуль
$\Gamma({\mathcal M},\End(E))$. Этот модуль можно рассматривать как
подалгебру в $\Diff(\Gamma({\mathcal M},E))$, состоящую из
"алгебраических", т.е. не содержащих дифференцирования, операторов.
Продолжим $\mu_E$ до отображения
$$\mu_E:\Gamma({\mathcal M},\End(E))\underset{{\mathfrak F}}{\otimes}
T({\mathcal M})\to\Diff(\Gamma({\mathcal M},E)),$$
потребовав, чтобы оно было перестановочно с алгебраическими
операторами (в очевидном смысле). После этого можно определить
геометрический символ обычным способом:
$$\sigma_E:\Diff(\Gamma({\mathcal M},E))\to
\Gamma({\mathcal M},\End(E))\underset{{\mathfrak F}}{\otimes}
S({\mathcal M}),\,\mu_E\circ\sigma_E=\id.$$
Оператор ковариантной симметризации
$$\Theta_E=\sigma_E\circ\mu_E:
\Gamma({\mathcal M},\End(E))\underset{{\mathfrak F}}{\otimes}
T({\mathcal M})\to
\Gamma({\mathcal M},\End(E))\underset{{\mathfrak F}}{\otimes}
S({\mathcal M})$$
также оказывается перестановочным с операторами из
$\Gamma({\mathcal M},\End(E))$, поэтому достаточно рассмотреть его ограничение
$$\Theta_E:T({\mathcal M})\to
\Gamma({\mathcal M},\End(E))\underset{{\mathfrak F}}{\otimes}
S({\mathcal M}).$$
Вычисления в этом случае существенно усложняются, поскольку $\Theta_E$
зависит не только от аффинной связности, но и от связности на расслоении $E$.

\section{Симметризация специальных тензоров}

В этом параграфе рассматривается алгоритм вычисления тензоров вида
$\Theta( v^{\otimes n}\otimes w\otimes v^{\otimes m})$ на многообразии
с аффинной связностью без кручения.
На практике симметризация осуществляется при помощи коммутационных
соотношений, таких как равенство
$$\nabla^3_{u,v,w}-\nabla^3_{v,u,w}=-\nabla_{R(u,v)w},$$
где $R$ - тензор кривизны.
Общие соотношения были получены автором в [6].
Нам, однако, потребуются соотношения достаточно специального вида.
Введем следующие обозначения:
$$q_2w=R(w,v)v,\,q_3w=(\nabla_v R)(w,v)v,\,q_4w=(\nabla^2_{v,v}R)(w,v)v;$$
в общем случае
$$q_n=q_n(v):w\mapsto (v^{\otimes n-2}\cdot\nabla^{n-2}R)(w,v)v,\,n\ge 2.$$
Мы считаем $q_n\in\End_{\Rs}({\mathcal V})$ операторами,
зависящими от $v\in{\mathcal V}$ как от параметра.

{\bf Лемма 4}
{\it Пусть $n,m\ge 0$, тогда для всяких полей $v,w\in{\mathcal V}$
на многообразии с симметричной связностью имеет место равенство
$$\mu(v^{\otimes n}\otimes w\otimes v^{\otimes m+1}-
v^{\otimes n+1}\otimes w\otimes v^{\otimes m}+$$
$$+\sum_{k=0}^{n}\sum_{l=1}^{m}\binom{n}{k}
v^{\otimes n+l-k-1}\otimes q_{k+2} w\otimes v^{\otimes m-l})=0,
\eqno{(7)}$$
где $q_k=q_k(v)$.
}

Например, при $n=0, m=1$
$$\mu(w\otimes v^2-v\otimes w\otimes v+q_2w)=
\nabla^3_{w,v,v}-\nabla^3_{v,w,v}+\nabla_{R(w,v)v}=0.$$
При $m=0$ сумма справа не имеет смысла, и мы условно считаем
её равной нулю:
$$\mu(v^{\otimes n}\otimes w\otimes v-
v^{\otimes n+1}\otimes w)=0,\,n\ge 0.$$

Доказательство леммы мы начнем со следующего равенства:
$$\mu(
w\otimes v\otimes u_1\otimes\dots\otimes u_m-
v\otimes w\otimes u_1\otimes\dots\otimes u_m+
\sum_{l=1}^{m}u_1\otimes\dots\otimes u_{l-1}\otimes
R(w,v)u_{l}\otimes u_{l+1}\otimes\dots\otimes u_m)=0,$$
где $u_1,\dots,u_m\in{\mathcal V}$.
Указанное равенство получается применением хорошо известного
правила перестановки индексов ко второй ковариантной производной
тензора $\nabla^m f,\,f\in{\mathfrak F}$.
При $u_1=\dots=u_m=v$ мы имеем
$$\mu(w\otimes v^{\otimes m+1}-
v\otimes w\otimes v^{\otimes m}+
\sum_{l=1}^{m}v^{\otimes l-1}\otimes q_2(v)w\otimes v^{\otimes m-l})=0.$$
Таким образом, мы получили утверждение леммы для $n=0$.

В общем случае утверждение можно получить, применяя к последнему
равенству  оператор $\mu(v^{\otimes n})$.
Для этого, однако, придется воспользоваться искусственным
приемом. Предположим вначале, что $\nabla_v v=\nabla_vw=0$. Это означает, что
интегральные кривые поля $v$ являются геодезическими, а поле $w$
ковариантно постоянно вдоль них. Согласно лемме 1 из [1],
$\mu_E(v^{\otimes k})=\nabla_v^k$ при $k>0$ для любого расслоения $E$.
Следовательно,
$$\mu_T(v^{\otimes k})w\otimes v^{\otimes m+1}=
\mu_T(v^{\otimes k})v\otimes w\otimes v^{\otimes m}=0,$$
$$\mu_T(v^{\otimes k})
(v^{\otimes l-1}\otimes q_2(v)w\otimes v^{\otimes m-l})=
v^{\otimes l-1}\otimes q_{k+2}(v)w\otimes v^{\otimes m-l}.$$
Эти равенства достаточно очевидны из геометрических соображений
(кроме того, их также можно получить как формальное следствие (4)).

Применяя формулу (4) к левой части равенства
$$\mu(v^{\otimes n})\mu(w\otimes v^{\otimes m+1}-
v\otimes w\otimes v^{\otimes m}+
\sum_{l=1}^{m}v^{\otimes l-1}\otimes q_2(v)w\otimes v^{\otimes m-l})=0,$$
мы получим
$$\mu(\sum_{k=0}^n\binom{n}{k}v^{\otimes n-k}\otimes\mu_T(v^{\otimes k})
[w\otimes v^{\otimes m+1}-
v\otimes w\otimes v^{\otimes m}+
\sum_{l=1}^{m}v^{\otimes l-1}\otimes q_2(v)w\otimes v^{\otimes m-l}
])=0,$$
т.е.
$$\mu(v^{\otimes n}\otimes w\otimes v^{\otimes m+1}-
v^{\otimes n+1}\otimes w\otimes v^{\otimes m}+
\sum_{k=0}^{n}\sum_{l=1}^{m}\binom{n}{k}
v^{\otimes n+l-k-1}\otimes q_{k+2} w\otimes v^{\otimes m-l})=0,$$
что и требовалось.

Чтобы избавиться от дополнительного условия, заметим, что оператор $\mu$
в левой части (7) зависит в каждой точке $x\in{\mathcal M}$
только от значения полей $v,w$ в этой же точке, т.е. от векторов
$v_x,w_x\in T_x{\mathcal M}$. Как легко убедиться, в
окрестности $U\subset{\mathcal M}$ любой точки $x$ можно построить поля
$\tilde{v},\tilde{w}\in{\mathcal V}(U)$, удовлетворяющий нужным условиям,
причем $\tilde{v}_x=v_x,\tilde{w}_x=w_x$. Тогда равенство (7) будет
справедливо в $U$ для оператора, полученного заменой $v,w$ на
$\tilde{v},\tilde{w}$. В точке $x$ этот оператор совпадает с исходным,
поэтому равенство (7) выполняется в точке $x$. Так как точка
может быть произвольной, лемма доказана.

Ковариантная симметризация тензора, зависящего от двух переменных
(векторных полей) и линейного по одной из них может быть произведена
путем многократного применения формулы (7). На каждом шаге исходный многочлен
представляется в виде суммы трех слагаемых, одно из которых лежит в
$S({\mathcal M})$, другое в  $\ker(\mu,T({\mathcal M}))$, а третье
играет роль остатка. Тензорная степень остатка меньше степени исходного
многочлена как минимум на 2, поэтому число шагов не превосходит половины
степени многочлена, подвергаемого симметризации.

Предлагаемый алгоритм лучше всего проиллюстрировать примером. Рассмотрим моном
$v^{\otimes 4}\otimes w,\,(v,w\in{\mathcal V})$.
Очевидно,
$$v^{\otimes 4}\otimes w=\frac{1}{5}(v^{\otimes 4}\otimes w+
v^{\otimes 3}\otimes w\otimes v+v^{\otimes 2}\otimes w\otimes v^{\otimes 2}+
v\otimes w\otimes v^{\otimes 3}+w\otimes v^{\otimes 4})-
\frac{1}{5}(w\otimes v^{\otimes 4}-v\otimes w\otimes v^{\otimes 3})-$$
$$-\frac{2}{5}(v\otimes w\otimes v^{\otimes 3}-
v^{\otimes 2}\otimes w\otimes v^{\otimes 2})-
\frac{3}{5}(v^{\otimes 2}\otimes w\otimes v^{\otimes 2}-
v^{\otimes 3}\otimes w\otimes v)-
\frac{4}{5}(v^{\otimes 3}\otimes w\otimes v-v^{\otimes 4}\otimes w).$$

Чтобы упростить обозначения, будем писать $A\equiv B$, если
$\mu(A-B)=0$.
(Подразумевая $A\equiv B \mod\ker(\mu,T({\mathcal M}))$.
Следует помнить, что поскольку отображение $\mu$ не является гомоморфизмом
алгебр, его ядро - не идеал).
Применяя лемму 4 при $n+m=4$, мы получим
$$w\otimes v^{\otimes 4}-v\otimes w\otimes v^{\otimes 3}+
(q_2w\otimes v^{\otimes 2}+v\otimes q_2w\otimes v+
v^{\otimes 2}\otimes q_2w)\equiv 0,$$
$$v\otimes w\otimes v^{\otimes 3}-
v^{\otimes 2}\otimes w\otimes v^{\otimes 2}+
(v\otimes q_2w\otimes v+v^{\otimes 2}\otimes q_2w+
q_3w\otimes v+v\otimes q_3w)\equiv 0,$$
$$v^{\otimes 2}\otimes w\otimes v^{\otimes 2}-
v^{\otimes 3}\otimes w\otimes v+(v^{\otimes 2}\otimes q_2w+
2v\otimes q_3w+q_4w)\equiv 0,$$
$$v^{\otimes 3}\otimes w\otimes v-v^{\otimes 4}\otimes w\equiv 0.$$

Таким образом,
$$v^{\otimes 4}\otimes w-S_1\equiv Q_1,$$
где
$$S_1=\frac{1}{5}(v^{\otimes 4}\otimes w+
v^{\otimes 3}\otimes w\otimes v+v^{\otimes 2}\otimes w\otimes v^{\otimes 2}+
v\otimes w\otimes v^{\otimes 3}+w\otimes v^{\otimes 4})+
\frac{3}{5}q_4w+$$
$$+\frac{1}{5}(q_2w\otimes v^{\otimes 2}+v\otimes q_2w\otimes v+
v^{\otimes 2}\otimes q_2w)+\frac{2}{5}(q_3w\otimes v+v\otimes q_3w)
\in S({\mathcal M}),$$
$$Q_1=v^{\otimes 2}\otimes q_2w+\frac{2}{5}v\otimes q_2w\otimes v+
\frac{6}{5}v\otimes q_3w.$$

Многочлен $Q_1$ играет роль остатка.
Заметим, что он является суммой мономов вида
$v^{\otimes i}\otimes w^{\prime}\otimes v^{\otimes j}$ (где
$w^{\prime}=q_dw,\,d=2,3$). Поэтому к нему может быть
применена та же процедура:

$$Q_1=
\frac{7}{15}(q_2w\otimes v^{\otimes 2}+v\otimes q_2w\otimes v+
v^{\otimes 2}\otimes q_2w)+$$
$$+\frac{3}{5}(q_3w\otimes v+v\otimes q_3w)-
\frac{7}{15}(q_2w\otimes v^{\otimes 2}-v\otimes q_2w\otimes v)-
\frac{8}{15}(v\otimes q_2w\otimes v-v^{\otimes 2}\otimes q_2w)-
\frac{3}{5}(q_3w\otimes v-v\otimes q_3w).$$
Далее,
$$q_2w\otimes v^{\otimes 2}-v\otimes q_2w\otimes v+q_2^2w\equiv 0,$$
$$v\otimes q_2w\otimes v-v^{\otimes 2}\otimes q_2w\equiv
q_3w\otimes v-v\otimes q_3w\equiv 0,$$
откуда
$$Q_1\equiv S_2,$$
где
$$S_2=\frac{7}{15}(q_2w\otimes v^{\otimes 2}+v\otimes q_2w\otimes v+
v^{\otimes 2}\otimes q_2w)+
\frac{3}{5}(q_3w\otimes v+v\otimes q_3w)+
\frac{7}{15}q_2^2w\in S({\mathcal M}).$$
Как мы видим, на втором шаге остаток равен нулю, и мы получаем ответ:
$$\Theta (v^{\otimes 4}\otimes w)=S_1+S_2.$$

\section{Двойное экспоненциальное отображение}

Следуя [2], определим двойное экспоненциальное отображение \newline
${\rm exp}_x:T_x{\mathcal M}\times T_x{\mathcal M}\to {\mathcal M}$
согласно правилу
$${\rm exp}_x(v,w)=
{\rm exp}_{{\rm exp}_x(v)}(I_{{\rm exp}_x(v)}^xw),$$
где символом $I_y^x:T_x{\mathcal M}\to T_y{\mathcal M}$
обозначен оператор параллельного переноса вдоль геодезической.
Вообще говоря, эта функция, как и обычное экспоненциальное отображение,
определена лишь в окрестности нуля. В случае
полного аффинного многообразия ${\mathcal M}$ её можно считать
всюду определённой, полагая, что вектор $w\in T_x{\mathcal M}$
переносится вдоль геодезической ${\rm exp}(x,tv),\,0\le t\le 1.$
Нас, однако, интересует лишь её локальное поведение.

Поскольку экспоненциальное отображение локально является диффеоморфизмом,
в окрестности нуля однозначно определена гладкая функция
$$h_x:T_x{\mathcal M}\times T_x{\mathcal M}\to T_x{\mathcal M},\quad
{\rm exp}_x(h_x(v,w))={\rm exp}_x(v,w),\,v,w\in T_x{\mathcal M}.\eqno{(8)}$$
С точки зрения геометрии, $h_x$ является функцией перехода от нормальных
координат в точке ${\rm exp}_x(v)$ к нормальным координатам в точке $x$
(при фиксированном $v$). При этом касательные пространства в разных
точках отождествляются с помощью параллельного переноса вдоль геодезической.
Можно также называть эту функцию двойным экспоненциальным отображением
(в нормальных координатах).

В отличие от ${\rm exp}_x$, функция $h_x$ является отображением
линейных пространств, поэтому можно ставить вопрос о вычислении ее ряда
Тейлора. Из геометрических соображений ясно, что слагаемые этого ряда
будут выражаться через тензоры кривизны и кручения вместе с их
ковариантными производными.

Далее мы полагаем, что связность симметрична, т.е. что её тензор кручения
равен нулю. Первые члены ряда Тейлора для этого случая были найдены автором
непосредственно из определения [2]:
$$h_x(v,w)=v+w+\frac{1}{6}R(w,v)v+\frac{1}{3}R(w,v)w+
\frac{1}{12}(\nabla_vR)(w,v)v+\frac{1}{24}(\nabla_wR)(w,v)v+$$
$$+\frac{5}{24}(\nabla_vR)(w,v)w+\frac{1}{12}(\nabla_wR)(w,v)w+
r(v,w),\eqno{(9)}$$
где $r(v,w)$ содержит однородные члены степени 5 и выше.
Ряд (9) напоминает известный ряд Кэмпбелла-Хаусдорфа с той разницей, что
выражается он не через скобку Ли, а через производные
тензора кривизны $\nabla^nR,\,n\ge 0$, которые формально
являются полилинейными $n+3$ - арными операциями.

Удобно считать $v$ и $w$ не векторами, а векторными полями,
рассматривая вместо $h_x$ функцию
$$h:{\mathcal V}\times{\mathcal V}\to{\mathcal V},\,
h(v,w)_x=h_x(v_x,w_x),\,x\in{\mathcal M}.$$
Конечно, поле $h(v,w)$ может оказаться не везде определённым.
Чтобы оно было определено, достаточно, чтобы поля $v,w$ имели компактный
носитель и не были слишком велики. Поскольку вектор в точке всегда
можно продолжить до гладкого поля с компактным носителем,
переход от $h_x$ к $h$ ничего не меняет по существу, но зато позволяет
использовать наработанную технику.

Нашей целью будет вычисление линейной по второй переменной части ряда
Тейлора функции $h$, или, что то же самое, функции $h_x$.
Определим зависящий от переменной $v\in{\mathcal V}$ линейный оператор
$$H(v):{\mathcal V}\to{\mathcal V},\,
H(v)w=\frac{d}{dt}h(v,tw)\bigg{|}_{t=0},\,v,w\in{\mathcal V}.$$
Оператор $H(v)$ может быть разложен в ряд Тейлора:
$$H(v)=\sum_{n=0}^{\infty}\frac{1}{n!}H_n(v),$$
где $H_n$ -- однородные слагаемые степени $n$, т.е.
$H_n(tv)=t^nH_n(v),\,t\in\Rs$.
Поскольку $h(0,w)=w,$ $H_0$ -- тождественный оператор. Из (9) следует, что
$$H_1(v)=0,\,H_2(v)=\frac{1}{3}q_2(v),\,H_3(v)=\frac{1}{2}q_3(v).$$
В принципе, все операторы $H_n$ могут быть найдены непосредственно
из определения двойного экспоненциального отображения. Но такой подход приводит
к чрезвычайно трудоемким вычислениям уже при сравнительно небольших $n$.
Как мы увидим, существует более простой путь.

Воспользуемся теоремой 1 [2], которую можно записать в виде
$$\mu(e^{\otimes v}\otimes e^{\otimes w})=\mu(e^{\otimes h}).$$
Здесь левая часть равенства понимается как ряд Тейлора правой части
по переменным $v,w$,
а знак $e^{\otimes v}=\sum_{n=0}\frac{1}{n!}v^{\otimes n}$ обозначает
формальную экспоненту в тензорной записи.
Поскольку, очевидно, $h^{\otimes n}\in S({\mathcal M})$ для любого $n$,
это означает, что формально
$$\Theta:e^{\otimes v}\otimes e^{\otimes w}\mapsto e^{\otimes h}$$
по определению оператора симметризации.
Если мы обозначим через $\pi_1$ естественный проектор
$S({\mathcal M})=\bigoplus_{n=0}^{\infty}S^n({\mathcal M})$
на $S^1({\mathcal M})={\mathcal V},$
то будем иметь равенство
$h=\pi_1\circ\Theta (e^{\otimes v}\otimes e^{\otimes w}).$
Сравнивая однородные по $v$ слагаемые в левой и правой частях,
мы получаем
$$H_n(v)w=\pi_1\circ\Theta( v^{\otimes n}\otimes w).$$
Таким образом, для вычисления $H_n$
требуется ковариантная симметризация монома
$v^{\otimes n}\otimes w.$ Она может быть проделана
при помощи алгоритма, описанного в предыдущем параграфе.
Например, при $n=4$ мы имеем
$$H_4(v)w=\pi_1\circ\Theta (v^{\otimes 4}\otimes w)=
\pi_1(S_1+S_2)=\frac{7}{15}q_2^2w+\frac{3}{5}q_4w,$$
так что
$$H_4(v)=\frac{7}{15}q_2^2+\frac{3}{5}q_4.$$

Операторы $H_n$ при $n>4$ могут быть найдены тем же методом,
но за большее число шагов. При этом симметричные слагаемые
степени выше первой, которые нас не интересуют,
на каждом шаге можно просто опускать.
Описанный алгоритм был реализован автором в виде программы,
что позволило ему найти $H_n$ для  $n\le 30$.
Выражение для $H_n$ содержит операторы вида $q_{d_1}\dots q_{d_m}$,
удовлетворяющие условию $d_1+\dots+d_m=n$.
Как можно убедиться,  число различных операторов этого вида
равно $\varphi_{n-1}$,
т.е. $n-1$ - му числу Фибоначчи
($\varphi_{1}=\varphi_{2}=1,\,\varphi_{n+1}=\varphi_{n}+\varphi_{n-1},n>1$).
Например, выражение для $H_{20}$ содержит $\varphi_{19}=4181$
слагаемых, так что приводить его здесь едва ли уместно.
Результаты вычислений для $5\le n\le 10$ приведены ниже

$$H_5=\frac{2}{3}q_5+q_3q_2+\frac{4}{3}q_2q_3,$$
$$H_6=\frac{5}{7}q_6+\frac{11}{7}q_4q_2+\frac{25}{7}q_3^2+
\frac{18}{7}q_2q_4+\frac{31}{21}q_2^3,$$
$$H_7=\frac{3}{4}q_7+\frac{13}{6}q_5q_2+\frac{27}{4}q_4q_3+
\frac{33}{4}q_3q_4+\frac{25}{6}q_2q_5+\frac{17}{4}q_3q_2^2+
\frac{31}{6}q_2q_3q_2+\frac{73}{12}q_2^2q_3,$$
$$H_8=\frac{160}{9}q_3^2q_2+\frac{140}{9}q_3q_5
+\frac{182}{9}q_3q_2q_3+\frac{43}{5}q_4q_2^2+\frac{106}{9}q_2q_4q_2+
\frac{226}{9}q_2q_3^2+$$
$$+\frac{239}{15}q_2^2q_4+
\frac{127}{15}q_2^4+\frac{55}{9}q_2q_6+\frac{7}{9}q_8+
\frac{25}{9}q_6q_2+\frac{98}{9}q_5q_3+\frac{91}{5}q_4^2,$$
$$H_9=\frac{4}{5}q_9+\frac{17}{5}q_7q_2+\frac{168}{5}q_5q_4+
\frac{196}{5}q_4q_5+\frac{209}{5}q_4q_3q_2+\frac{232}{5}q_4q_2q_3+
\frac{74}{5}q_5q_2^2+$$
$$+47q_3q_4q_2+26q_3q_6+60q_3q_2q_4+31q_3q_2^3
+98q_3^3+\frac{168}{5}q_2^2q_5+\frac{226}{5}q_2^3q_3+$$
$$+\frac{42}{5}q_2q_7+
\frac{184}{5}q_2q_3q_2^2+\frac{378}{5}q_2q_3q_4+66q_2q_4q_3+22q_2q_5q_2+
\frac{197}{5}q_2^2q_3q_2+16q_6q_3,$$
$$H_{10}=\frac{9}{11}q_{10}+\frac{982}{11}q_5q_2q_3+\frac{2702}{33}q_5q_3q_2+
\frac{763}{33}q_6q_2^2+\frac{896}{11}q_5^2+\frac{612}{11}q_6q_4+\frac{243}{11}q_7q_3+$$
$$+\frac{133}{33}q_8q_2+\frac{4610}{33}q_2q_5q_3+
\frac{1205}{33}q_2q_6q_2+\frac{3358}{33}q_2q_4q_2^2+
\frac{2472}{11}q_2q_4^2+\frac{364}{33}q_2q_8+$$
$$+\frac{441}{11}q_3q_7+\frac{810}{11}q_4q_6+
\frac{7262}{33}(q_2q_3)^2+\frac{6494}{33}q_2q_3^2q_2+\frac{5936}{33}q_2q_3q_5+
\frac{3240}{11}q_3q_4q_3+\frac{3290}{33}q_3q_5q_2+$$
$$+\frac{3627}{11}q_3^2q_4+
\frac{5173}{33}q_3^2q_2^2+\frac{1554}{11}q_3q_2q_5+\frac{5348}{33}(q_3q_2)^2+
\frac{2028}{11}q_3q_2^2q_3+\frac{7945}{33}q_2^2q_3^2+$$
$$+\frac{2835}{11}q_4q_3^2+
\frac{2050}{33}q_2^2q_6+\frac{3787}{33}q_2^2q_4q_2+\frac{1383}{11}q_4^2q_2+
\frac{1636}{11}q_2^3q_4+\frac{1692}{11}q_4q_2q_4+\frac{2555}{33}q_2^5+
\frac{855}{11}q_4q_2^3.$$

В работе [5], посвящённой анализу на многообразии со связностью,
Шарафутдинов ввел многочлены $\rho^{\alpha,\beta}$, играющие роль своего
рода структурных констант в алгебре дифференциальных операторов.
Как было показано автором в [2], между этими многочленами
и двойным экспоненциальным отображением имеется связь,
которую можно описать при помощи производящей функции:
$$\sum_{\alpha,\beta}i^{|\alpha|+|\beta|}
\frac{v^{\alpha}w^{\beta}}{\alpha!\beta!}\rho^{\alpha,\beta}(x,\xi)=
e^{i\langle h(v,w)-v-w,\xi\rangle}.\eqno{(10)}$$

Сравнивая коэффициенты при $v^{\alpha}w^{\beta}$ для $|\beta|=1$,
мы получим равенство
$$\sum_{|\alpha|=n}\sum_{k}
\frac{v^{\alpha}}{\alpha!}w^{k}
\rho^{\alpha,\langle k\rangle}(x,\xi)=
\frac{i^{-n}}{n!}\langle H_n(v)w,\xi\rangle$$
или, в индексных обозначениях,
$v^{j_1}\dots v^{j_n}w^k
\rho^{\langle j_1\dots j_n\rangle,\langle k\rangle}(x,\xi)=
i^{-n}\langle H_n(v)w,\xi\rangle.$
Таким образом, многочлены  $\rho^{\alpha,\beta}$ при $|\beta|=1$
непосредственно связаны с оператором  $H_{|\alpha|}$.
Например, для $n=3$ мы  имеем
$$v^jv^kv^lw^m
\rho^{\langle jkl\rangle,\langle m\rangle}(x,\xi)=
i\langle H_3(v)w,\xi\rangle.$$
Так как $\langle H_3(v)w,\xi\rangle=\langle\frac{1}{2}q_3(v)w,\xi\rangle=
\frac{1}{2}v^jv^kv^lw^m\nabla_jR_{kml}^p\xi_p,$
мы, в обозначениях [5], получаем равенство
$$\rho^{\langle jkl\rangle,\langle m\rangle}(x,\xi)=-
\frac{i}{2}\sigma(jkl)\nabla_jR_{klm}^p\xi_p.$$

Вычисление $\rho^{\alpha,\beta}$ для $|\beta|>1$ при помощи
производящей функции сводится к нахождению слагаемых ряда Тейлора
функции $h$, порядок которых по второй переменной не превосходит $|\beta|$.
В принципе, эта задача может быть решена тем же методом,
с использованием вместо (7) общих коммутационных соотношений [6].
Однако на практике соответствующие вычисления существенно сложнее.


Гаврилов Алексей Владимирович,\\
E-mail:gavrilov19@gmail.com\\

\end{document}